\documentclass[12pt]{amsart}
\usepackage{amscd}
\usepackage{amssymb}
\usepackage{hyperref}

\def\depth{\operatorname{depth}}
\def\hom{\operatorname{Hom}}

\def\reg{\operatorname{reg}}
\def\pd{\operatorname{pd}}
\def\max{\operatorname{max}}
\def\reltype{\operatorname{reltype}}
\def\grade{\operatorname{grade}}

\newcommand{\m}{\mathfrak m}

\newcommand{\ZZ}{\mathbb Z}
\newcommand{\NN}{\mathbb N}

\newcommand{\F}{\mathcal F}
\newcommand{\limk}{\underset{\underset{k}{\longrightarrow}}{\lim~}}

\newtheorem{lemma}{Lemma}[section]
\newtheorem{corollary}[lemma]{Corollary}
\newtheorem{theorem}[lemma]{Theorem}
\newtheorem{proposition}[lemma]{Proposition}

\newtheorem{remark}[lemma]{Remark}
\newtheorem{example}[lemma]{Example}
\newtheorem{question}[lemma]{Question}

\advance\headheight1.15pt
\begin{document}
\title[Regularity and Gorensteinness of fiber cone]
{Castelnuovo-Mumford regularity and Gorensteinness of fiber cone}

\author{A. V. Jayanthan}
\address{Department of Mathematics, Indian Institute of Technology
Madras, Chennai, INDIA - 600036.}
\thanks{$^\dag$ Supported by the Council of Scientific and Industrial
Research, India}
\thanks{AMS Classification 2000: 13H10, 13C15}
\email{jayanav@iitm.ac.in}
\author[Ramakrishna Nanduri]{Ramakrishna Nanduri$^\dag$}
\address{Department of Mathematics, Indian Institute of Technology
Madras, Chennai, INDIA - 600036.}
\email{nandurirk@gmail.com}

\begin{abstract}
In this article, we study the Castelnuovo-Mumford regularity and
Gorenstein properties of the fiber cone. We obtain upper bounds for
the Castelnuovo-Mumford regularity of the fiber cone and obtain
sufficient conditions for the regularity of the fiber cone to be equal
to that of the Rees algebra. We obtain a formula for the canonical
module of the fiber cone and use it to study the Gorenstein property
of the fiber cone.
\end{abstract}

\maketitle
\section{Introduction}
Let $(A,\m)$ be a commutative Noetherian local ring and $I$ be an
ideal of $A$.  The graded algebras, the Rees algebra, $R(I) :=
\oplus_{n\geq0}I^nt^n \subset R[t]$, the associated graded ring, $G(I)
:= \oplus_{n\geq0}I^n/I^{n+1} \cong R(I)/IR(I)$ and the fiber cone,
$F(I) := \oplus_{n\geq0}I^n/\m I^n$ are together known as the blowup
algebras associated to $I$.
In this article, our aim is to study the Castelnuovo-Mumford
regularity and the Gorensteinness of the fiber cone. We do this by
relating them with the corresponding properties of certain other
graded modules.  Cortadellas and Zarzuela, in a series papers, used
certain graded modules associated to filtrations of modules to study
the depth properties of the fiber cone \cite{cz1}, \cite{cz2},
\cite{co2}. We use these graded modules and other blowup algebras to
study the regularity and the Gorensteinness of the fiber cone.

For a standard graded algebra $S = \oplus_{n\geq0}S_n$ over a
commutative Noetherian ring $S_0$ and a finitely generated graded
$S$-module $M = \oplus_{n\geq0}M_n$, define
$$
a(M) := \left\{
\begin{array}{ll}
\max\{n \; | \; M_n \neq 0\} & \mbox{ if } M \neq 0 \\
-\infty & \mbox{ if } M = 0.
\end{array}
\right.
$$
For $i \geq 0$, set
$$
a_i(M) := a(H^i_{S_+}(M)),
$$
where $S_+$ denotes the ideal of $S$ generated by the homogeneous
elements of positive degree and $H^i_{S_+}(M)$ denotes the $i$-th
local cohomology module of $M$ with respect to the ideal $S_+$. The
{\em Castelnuovo-Mumford regularity} (or {\em regularity}) of
$M$ is defined as the number
$$
\reg(M) := \max\{a_i(M) + i \; | \; i \geq 0 \}.
$$
Let $(A,\m)$ be a local ring and $I$ be any ideal.
The Castelnuovo-Mumford regularity of $R(I)$ and $G(I)$ have been well
studied in the past. Ooishi proved that $\reg R(I)  = \reg G(I)$,
\cite{o} (see also \cite{t1}). In \cite{t1}, Trung studied the
vanishing behavior of the local cohomology modules of the
associated graded ring and the Rees algebra and derived that for any
ideal in a Noetherian local ring $\reg R(I) = \reg G(I)$. It
can easily be seen that such an equality is not true in the case of
the fiber cone and the associated graded ring (see Section 2). In
Section 2, we prove that for any ideal of analytic spread one in
a Noetherian local ring, the regularity of the fiber cone is bounded
above by the regularity of the associated graded ring. If the ideal
contains a regular element, we show that the equality holds in the
above case (Theorem \ref{thm7}). We also prove that,
under some assumptions, the regularity of the fiber cone is
bounded below by the regularity of the associated graded ring and
obtain certain sufficient conditions for the equality.

In Section 3, we study the Gorenstein property of the fiber cone.
The Gorenstein property of the Rees algebra and the
associated graded ring has been very well studied, see for example
\cite{gn}, \cite{hrs}, \cite{hku}, \cite{eh}, \cite{tvz}.  It is known
that the Gorenstein fiber cones behave differently. For
example, it is known that, unlike in
the case of associated graded ring, fiber cone can be Gorenstein
without the ambient ring being Gorenstein. Also, it can easily be seen
that the symmetry of the Hilbert series of the fiber cone does not
assure the Gorensteinness for fiber cone (\cite[Example 6.2]{jpv}). In
Section 3, we obtain an expression for the canonical module of the
fiber cone (Proposition \ref{pro4}). Using the structure of the
canonical module, we obtain a necessary and sufficient condition for
the Gorenstein property of the fiber cone (Theorem \ref{gor-char}). We
also obtain an upper bound for the regularity of the canonical module
of $F(I)$. We end the article by proving that the multiplicity of the
canonical module of the fiber cone is strictly less than the multiplicity
of the canonical module of the associated graded ring, except for the
maximal ideal.  Throughout this article, $(A,\m)$ will always denote a
Noetherian local ring of dimension $d$ and with infinite residue
field. All the computations in this article have been performed using
the computer algebra package CoCoA, \cite{co}. \\
{\bf Acknowledgement :} We would like to sincerely thank the referee for
a thorough reading of the manuscript, pointing out some mistakes and 
suggesting some improvements.

\section{Castelnuovo-Mumford regularity of the fiber cone}

In this section we study the regularity of the fiber cone. Unlike in
the case of the associated graded ring and the Rees algebra, there is
no equality between the regularities of the fiber cone and the Rees
algebra. For example if $A$ is a non-Buchsbaum ring and $I$ is
generated by a system of parameters and not a $d$-sequence (such a
sequence exists by \cite[Proposition 1.7]{h}), then from Corollary 5.2
in \cite{t1} we have $\reg R(I) > 0$. Since $I$ is generated by a
system of parameters $\reg F(I) = 0$. Therefore $\reg F(I) < \reg
R(I)$. We show that under some hypothesis, the regularity of the fiber
cone is at most that of the Rees algebra. We also provide some
sufficient conditions for the equality.

Let $A$ be a Noetherian local ring of dimension $d > 0$ and
$I \subset A$ be an ideal. Consider the filtration
$$
\F : A \supset \m \supset \m I \supset \m I^2 \supset \cdots
$$
Let $R(\F) := A \oplus \m t \oplus \m It^2 \oplus \m I^2t^3 \oplus
\cdots$. Then $R(\F)$ is a finitely generated graded $R(I)$-module.
Consider the exact sequences of $R(I)$-modules:
\begin{equation}\label{eqn3}
 0 \rightarrow R(I) \rightarrow R(\mathcal{F}) \rightarrow \m G(I)(-1) \rightarrow 0
\end{equation}
\begin{equation} \label{eqn6}
 0 \rightarrow \m G(I) \rightarrow G(I) \rightarrow F(I) \rightarrow 0
\end{equation}

We use the above two exact sequences and the corresponding long exact
sequence of local cohomology modules to the study the vanishing
properties of the local cohomology modules of the fiber cone.
Throughout this section we assume that $I$ is an ideal of analytic
spread := $\dim F(I) = \ell > 0$. We begin this section with a remark on
the top local cohomology modules of $R(I)$ and $R(\F)$.

\begin{remark} \label{rmk2}
Suppose $\underline{x}=x_1,\ldots,x_\ell$ generates a minimal
reduction of $I$. Note that $R(I)_+$ is generated radically by $\ell$
elements. Denote by $\underline{x^k}$ the sequence
$x_1^k,\ldots,x_\ell^k$.  Then for all $n \in \ZZ$,
$$\left [H^\ell_{R(I)_+}(R(I))
\right ]_n\cong \limk \frac{I^{\ell k+n}}{(\underline{x}^k)
I^{(\ell-1)k+n}} \mbox{ and }
\left [ H^\ell_{R(I)_+}(R(\mathcal{F}))\right ]_n \cong \limk \frac{\m
I^{\ell k+n-1}}{(\underline{x}^k) \m I^{(\ell-1)k+n-1}}.$$
This implies that $a_\ell(R(\mathcal{F}))-1 \leq a_\ell(R(I))$.
\end{remark}

In the following theorem, we show that the regularity of the fiber
cone is bounded above by the regularity of the associated graded ring
when $\ell = 1$. For convenience, we will denote $H^i_{R(I)_+}(M)$ by
$H^i(M)$ for the rest of the section.
\begin{theorem} \label{thm7}
Let $(A, \m)$ be a Noetherian local ring and $I$ be an ideal of $A$
with $\ell= 1$. Then $\reg F(I) \leq \reg G(I)$. Furthermore, if
$\grade I = 1$, then $\reg F(I)=\reg G(I) = r(I),$ where $r(I)$
denotes the reduction number of $I$.
\end{theorem}
\begin{proof}
Since by hypothesis $\ell=1$, $I$ is generated by a single
element up to radical. Therefore $H^i(M)= (0)$ for all $i \geq 2$ and
for any finitely generated graded $R(I)$-module $M$. From the exact
sequence \eqref{eqn6}, we have the long exact sequence:
\begin{eqnarray*}
0 &\rightarrow& H^0(\m G(I)) \rightarrow H^0(G(I)) \rightarrow H^0(F(I)) \\
  &\rightarrow& H^1(\m G(I)) \rightarrow H^1(G(I)) \rightarrow H^1(F(I)) \rightarrow 0.
\end{eqnarray*}
It follows that $a_0(\m G(I)) \leq a_0(G(I))$ and $a_1(F(I)) \leq
a_1(G(I))$. From the exact sequence \eqref{eqn3}, we have the long
exact sequence:
\begin{eqnarray*}
0 &\rightarrow& H^0(R(I)) \rightarrow H^0(R(\mathcal{F})) \rightarrow H^0(\m G(I))(-1) \\
  &\rightarrow& H^1(R(I)) \rightarrow H^1(R(\mathcal{F})) \rightarrow
  H^1(\m G(I))(-1) \rightarrow 0 .
\end{eqnarray*}
Therefore $a_1(\m G(I)(-1)) = a_1(\m G(I))+1 \leq
a_1(R(\mathcal{F}))$.  From Remark \ref{rmk2}, it follows that
$a_1(\m G(I)) \leq a_1(R(I))$.  Since $G_+$ is radically generated
by one element, $H^1(G(I)) \neq 0$ and hence by Theorem 3.1 of
\cite{t1}, we get $a_1(R(I)) \leq a_1(G(I))$ so that $a_1(\m G(I))
\leq a_1(G(I))$. Therefore $\reg \m G(I) \leq \reg G(I)$. Now the
regularity behaviour under the exact sequence \eqref{eqn6} yields
that $$\reg F(I) \leq \max\{\reg G(I),\reg \m G(I)-1\}= \reg G(I).$$
\vskip 0.2cm Now assume that $\grade I = 1$. Then by (\cite{cz3},
page 764) we have $\reg F(I)=r_J(I)$, for any minimal reduction $J$
of $I$. Also we have $\reg G(I)= r_J(I)$ (see for example
(\cite{hz}, Proposition 3.6) ). Therefore $\reg F(I)=\reg G(I) =
r(I)$.
\end{proof}

Now we give a lower bound for the regularity of fiber cone under some
assumptions. In \cite{cz3}, Cortadellas and Zarzuela proved that if
the depth of the fiber cone and the associated graded ring is at least
$\ell -1$, then the regularities of these two algebras are equal. We
generalize this result in the following theorem and retrieve their
result in the above mentioned case.
For $x \in I \; \backslash \m I$, let $x^*$ denote
the image of $x$ in $I/I^2$ and $x^o$ denote the image of $x$ in $I/\m
I$.
\begin{theorem} \label{pro8}
Let $(A,\m)$ be a Noetherian local ring and $I$ be an ideal of $A$.
Suppose $\grade I = \ell$ and $\grade G(I)_+ \geq \ell-1$. Then
$\reg F(I) \geq \reg G(I)$.  Furthermore, if $\depth F(I) \geq
\ell-1$, then $\reg F(I)= \reg G(I)$.
\end{theorem}
\begin{proof} If $\ell=1$, then the proposition follows from Theorem
\ref{thm7}.  Suppose $\ell \geq 2$. Let $x_1,\ldots,x_\ell$ be a
minimal generating set for a minimal reduction $J$ of $I$ such that
$x_1^*,\ldots,x_\ell^* \in I/I^2$ is a filter regular sequence for
$G(I)$ and $x_1^o,\ldots,x_\ell^o \in I/\m I$ is a filter regular
sequence for $F(I)$. Since $\grade G(I)_+ \geq \ell -1, \;
x_1^*,\ldots,x_{\ell-1}^*$ is $G(I)$-regular.  Let ``-'' denote
modulo $(x_1,\ldots,x_{\ell-1})$. Then $G(\bar{I}) \cong
G(I)/(x_1^*,\ldots,x_{\ell-1}^*), \; F(\bar{I}) \cong
F(I)/(x_1^o,\ldots,x_{\ell-1}^o)$ and $\reg G(\bar{I})= \reg G(I)$.
Since $\dim F(\bar{I}) = 1$, by Theorem \ref{thm7} we get $\reg
F(\bar{I})= \reg G(\bar{I}) = \reg G(I)$. From \cite[Proposition
1.2]{ch}, it follows that $\reg F(I) \geq \reg F(\bar{I})$. This
implies that $\reg F(I) \geq \reg G(I)$. \vskip 2mm Now assume
$\depth F(I) \geq \ell-1$. Then $x_1^o,\ldots,x_{\ell-1}^o$ is
$F(I)$-regular. Then $\reg F(I)= \reg F(\bar{I})$. Therefore $\reg
F(I)= \reg G(I)$ as required.
\end{proof}

Now we give certain instances where the regularity of the fiber cone
is equal to the regularity of the Rees algebra or the associated
graded ring.
\begin{proposition}\label{pro6}
Let $(A,\m)$ be a Noetherian local ring and $I$ be an
ideal of $A$. If $\reg R(\mathcal{F}) \leq \reg R(I)$, then $\reg
F(I)= \reg R(I)$.
\end{proposition}
\begin{proof}
From the exact sequence \eqref{eqn3} and the fact that $\reg \m
G(I)(-1)=\reg \m G(I)+1$, it follows that $\reg \m G(I)+1 \leq \max\{
\reg R(I)-1, \reg R(\mathcal{F})\}$. Since $\reg R(\mathcal{F}) \leq
\reg R(I)$, the above inequality implies that $\reg \m G(I)+1 \leq \reg
R(I)$.  From the exact sequence \eqref{eqn6}, we get $\reg F(I) \leq
\max \{\reg \m G(I)-1, \reg G(I)\}$. Since $\reg R(I) = \reg G(I)$,
the above inequality implies that $\reg F(I) \leq \reg R(I)$.
\vskip 2mm
Now from the exact sequence \eqref{eqn6}, $\reg G(I)\leq \max \{\reg
\m G(I),\reg F(I)\}$. Since $\reg G(I)=\reg R(I)$ and $\reg \m G(I)
\leq \reg R(I)-1,$
the above inequality yields that
$$\max \{\reg \m G(I),\reg F(I)\}= \reg F(I) \mbox{ and }
\reg R(I) = \reg G(I) \leq \reg F(I).$$ Therefore we have
$\reg R(I)=\reg F(I)$. This completes the proof.
\end{proof}
Note that, if $\grade_{R(I)_+} R(\F) \geq \ell$, then
$H^i(R(\F)) = 0$ for $i < \ell$ and from Remark \ref{rmk2}, it follows
that $\reg R(\F) \leq \reg R(I)$. The next proposition gives yet
another instance of the equality of the regularity of these graded
algebras.
\begin{proposition}\label{pro9}
Let $(A,\m)$ be a Noetherian local ring and $I$ be an $\m$-primary
ideal of $A$ such that $\grade I > 0$. Suppose $I^{n_{0}}= \m
I^{n_{0}-1}$ for some $n_0 \in \NN$. Then $\reg F(I)= \reg G(I)$.
\end{proposition}
\begin{proof}
Since $I^{n_{0}}= \m I^{n_{0}-1}$ for some $n_0 \in
\NN$, it follows that $\m G(I)$ is Artinian. If $I = \m$, then the
assertion of the theorem follows trivially. If $I \neq \m$, then $\m
G(I) \neq 0$. Therefore $H^0(\m G(I))= \m G(I) \neq 0$
and $H^i(\m G(I))= 0$ for all $i > 0$.
From the exact sequence \eqref{eqn6} we have
$$0 \rightarrow H^0(\m G(I)) \rightarrow H^0(G(I)) \rightarrow
H^0(F(I)) \rightarrow 0$$
and $H^i(G(I)) \cong H^i(F(I))$ for all $i > 0$. Therefore $a_0(F(I))
\leq a_0(G(I))$ and $a_i(G(I))= a_i(F(I))$ for $i > 0$. If $\depth
G(I) >0$, then $H^0(G(I))=(0)$ and hence from the above exact sequence
we get $H^0(\m G(I))= 0$, which is a contradiction. Therefore $\depth
G(I)=0$.  Then by the Proposition 6.1 in \cite{t1}, we have
$$a_0(F(I)) \leq a_0(G(I)) < a_1(G(I)) = a_1(F(I)).$$
Therefore $\reg
F(I)= \max\{a_i(F(I))+i : i \geq 1 \}= \max\{a_i(G(I))+i : i \geq 1\}=
\reg G(I)$ as required.
\end{proof}

\begin{proposition}
Let $(A,\m)$  be a Noetherian local ring and $I$ be an ideal of $A$
such that $\grade I > 0$. Assume that $\m G(I)$ is a Cohen-Macaulay
$R(I)$-module of dimension $\ell$. Then
\begin{enumerate}
\item[(i)] $\reg F(I) \leq \reg R(I)$;
\item[(ii)] if $a_\ell(R(\mathcal{F}))-1 < a_\ell(R(I))$, then $\reg
F(I) = \reg R(I)$;
\item[(iii)] if $a_\ell(R(\mathcal{F}))-1 = a_\ell(R(I)),$  then $\reg
\m G(I) \leq \reg R(I)$ and $\reg F(I) \leq \reg R(I)$. Furthermore,
if $\reg \m G(I) < \reg G(I)$, then $\reg F(I) = \reg R(I)$.
\end{enumerate}
\end{proposition}
\begin{proof}
From the short exact sequence \eqref{eqn3}, there is a long exact
sequence of the local cohomology modules:
  $$\cdots \rightarrow H^i(R(I)) \rightarrow H^i(R(\mathcal{F})) \rightarrow H^i(\m G(I)(-1)) \rightarrow H^{i+1}(R(I)) \rightarrow \cdots.$$
Since $\m G(I)$ is Cohen-Macaulay, $H^i(\m G(I))=0$ for $i \neq \ell$.
Therefore, it follows from the above long exact sequence that
$H^i(R(I)) \cong H^i(R(\mathcal{F}))$ for $i < \ell$. This implies
that $a_i(R(\mathcal{F})) = a_i(R(I))$ for $i < \ell$. Using Remark \ref{rmk2},
we have $a_\ell(R(\mathcal{F}))-1 \leq a_\ell(R(I))$. Therefore (i)
follows from (ii) and (iii).
\vskip 0.2cm
If $a_\ell(R(\mathcal{F}))-1 < a_\ell(R(I))$, then
$a_\ell(R(\mathcal{F})) \leq a_\ell(R(I))$. Therefore $a_i(R(\mathcal{F}))
\leq a_i(R(I))$ for all $i$.  This implies that $\reg R(\mathcal{F})
\leq \reg R(I)$. Therefore from the Proposition \ref{pro6} we have
$\reg F(I)= \reg R(I)$. This proves (ii).
\vskip 0.2cm
Suppose $a_\ell(R(\mathcal{F}))-1 = a_\ell(R(I))$. Since
$a_i(R(\mathcal{F})) = a_i(R(I))$ for $i < \ell,$ we have
$\reg R(\mathcal{F}) \leq \reg R(I)+1$. Therefore from the  exact sequence
\eqref{eqn3}, it follows that
$$\reg(\m G(I)(-1)) = \reg \m G(I)+1 \leq \max\{ \reg R(I)-1, \reg
R(\mathcal{F})\} \leq \reg R(I)+1.$$
From the exact sequence \eqref{eqn6}, it follows that $\reg F(I) \leq \max \{\reg \m G(I)-1, \reg G(I)\} \leq  \reg R(I)$.
Thus $\reg F(I) \leq \reg R(I)$.
\vskip 0.2cm
Now assume that $\reg \m G(I) < \reg G(I)$. From the exact sequence
\eqref{eqn6} it follows that $\reg G(I)\leq \max \{\reg \m G(I),\reg
F(I)\}$.  Since $\reg \m G(I) < \reg G(I)$, the above inequality gives
that $\max \{\reg \m G(I),\reg F(I)\}= \reg F(I)$. Therefore $\reg
G(I) \leq \reg F(I)$.  Since $\reg G(I)= \reg R(I)$, we have $\reg
R(I) \leq \reg F(I)$. The other inequality is already proved.
Therefore $\reg F(I) = \reg R(I)$. This proves (iii).
\end{proof}

We conclude this section by giving some examples to illustrate the regularity behavior of the
fiber cone. The following example shows that the regularity of the
fiber cone can be strictly less than the regularity of the associated
graded ring even when $\reg F(I) > 0$. Let $R = \oplus_{n\geq 0}R_n$
be a finitely generated standard graded algebra. Then $R \cong
R_0[X_1, \ldots, X_m]/J$ for some $m$ and a homogeneous ideal $J$,
where $X_1, \ldots, X_m$ are indeterminates over $R_0$. Then the
relation type of $R$, denoted by $\reltype(R)$ is defined to be the
maximum degree of a minimal generating set of $J$.  It is known that
$\reltype(R) \leq \reg R + 1$, \cite{t1}. Let $k$ denote a field.
\begin{example}\label{reg-less}
Let $A = k[\![X,Y,Z]\!]/(X^2,Y^2,XYZ^2)$ and $\m = (x,y,z)$, where
$x=\bar{X}, y=\bar{Y}, z=\bar{Z}$ and $k$ is a field. Then $A$ is a one dimensional
Noetherian local ring. Let $I =(y,z)A$. Then $F(I) \cong k[Y,Z]/(Y^2)$. Therefore $\reg F(I)=1$. Let $\psi : A/I[U,V]\rightarrow G(I)$
be $A/I$-algebra homomorphism defined by $\psi(U)= yt$ and $\psi(V)= zt$. Then $ker(\psi)$ has a generator $(x+I)UV^2$. Therefore
$\reltype G(I) \geq 3$. Since $\reltype G(I) \leq \reg G(I)+1$ we have $\reg G(I) \geq 2$. Thus $\reg F(I) < \reg G(I)$. Here $\grade I = 0$.
\end{example}
\begin{example}
Let $A= k[\![x,y,z]\!]/J$, where $$J = (xy^3, xy^2z, xyz^2, xz^3, x^3z^2,
x^4, y^3z, x^3y, x^2y^2, y^4)$$ and $k$ is a field.  Let $I=
(\bar{x}^2,\bar{y},\bar{z})$ and $\m= (\bar{x},\bar{y},\bar{z})$. Then
$A$ is one dimensional non-Cohen-Macaulay Noetherian local ring and
$I$ is an $\m$-primary ideal satisfying $I^4= \m I^3$. Let $S=
k[U,V,W]$. Then $F(I) \cong \frac{S}{(U^2,UV^2,V^3W,UVW^2,UW^3,V^4)}$.
By using CoCoA, \cite{co}, one can see that the minimal free resolution
of $F(I)$ as an $S$-module is 
$$0 \rightarrow S(-6)^3 \rightarrow S(-4) \oplus S(-5)^7 \rightarrow
S(-2)\oplus S(-3) \oplus S(-4)^4 \rightarrow S \rightarrow 0.$$ 
From this the regularity of $F(I)$ is $3$. Therefore by the
Proposition \ref{pro9} we have $\reg F(I)= \reg G(I)=3$.  Note that
$\grade I=0$.
\end{example}
The following example shows that the reduction number of $I$ can be
strictly smaller than the regularity of $F(I)$.
\begin{example}
Let $A= k[\![X,Y,Z]\!]/J$, where $J=(X^4,XY^2Z,XYZ^2,YZ^4,Z^5)$. Let
$I=(x^3,y^2,z^2)$ and $\m=(x,y,z)$, where $x=\bar{X},y=\bar{Y}$ and
$z=\bar{Z}$.  Then $A$ is a non-Cohen-Macaulay Noetherian local ring.
Let $J =(y^2)$, then $I^4 = JI^3$. Therefore $r(I) \leq 3$. Let
$S=k[U,V,W]$.  Then $F(I)\cong \frac{S}{(U^2,UVW,W^3,V^2W^2,V^4W)}$.
By using CoCoA, \cite{co}, one can see that the minimal free
resolution of $F(I)$ as $S$-module is
\begin{eqnarray*}
0 &\longrightarrow& R(-6)^2 \oplus R(-7) \longrightarrow R(-4) \oplus
R(-5)^4 \oplus R(-6)^2 \\
&\longrightarrow& R(-2) \oplus R(-3)^2 \oplus R(-4) \oplus R(-5)
\longrightarrow R \longrightarrow 0.
\end{eqnarray*}
Hence $\reg(F(I)) = 4$. Therefore $r(I)  < \reg F(I)$.
\end{example}

Though we have proved that under certain conditions, $\reg F(I)$ is
bounded below by $\reg G(I)$, we have not been able to find an example
when this inequality is strict. Therefore we ask:

\begin{question}
Let $(A,\m)$ be a Noetherian local ring with infinite residue field
and $I$ an ideal such that $\grade I > 0$. Is $\reg F(I) \leq \reg
G(I)$?
\end{question}
\section{Gorenstein fiber cones}
In this section, we study the Gorenstein property of the fiber cone.
We begin by obtaining an expression for the canonical module of the
fiber cone. We fix the notation for this section. Throughout this
section, we assume that $G(I)$ and $F(I)$ are Cohen-Macaulay.
Let $\omega_{G(I)}$ and $\omega_{F(I)}$ denote the canonical modules of the
associated graded ring $G(I)$ and the fiber cone $F(I)$ respectively. For the
definition and basic properties of canonical modules, see \cite{bh}.
In the original manuscript, the result given below was proved in a
weaker form.  We would like to thank the referee for suggesting the
following improved form. 

\begin{proposition} \label{pro4}
Let $(A,\m)$ be a Noetherian local ring and $I$ be an $\m$-primary ideal
such that the associated graded ring $G(I)$ is Cohen-Macaulay. Let
$\omega_{G(I)}= \oplus_{n \in \ZZ} \omega_n$ and $\omega_{F(I)}$
be the canonical modules of $G(I)$ and $F(I)$ respectively. Then
\begin{enumerate}
\item[(1)] $\omega_{F(I)}\cong \oplus_{n \in \ZZ}(0:_{\omega_n}\m )$;
\item[(2)] $a(F(I))=a(G(I))= r-d$, where $r$ is the reduction number
  of $I$ with respect to any minimal reduction $J$ of $I$;
\item[(3)] for any $k \in \NN$, $a(F(I^k))= [\frac{a(F(I))}{k}]= [\frac{r-d}{k}]$;
\item[(4)] if $G(I)$ is Gorenstein, then
\end{enumerate}
$$\omega_{F(I)} \cong \bigoplus_{n \in \ZZ} \frac{(I^{n+r-d+1}:\m)\cap
I^{n+r-d}}{I^{n+r-d+1}}.$$
\end{proposition}
\begin{proof}
(1) Since $G(I)$ is Cohen-Macaulay and $F(I)=G(I)/\m G(I)$ is such
that $\dim G(I)=\dim F(I)= d$ we have by (\cite{hio}, Corollary (36.14))
that:
$$\omega_{F(I)} \cong Hom_{G(I)}(F(I), \omega_{G(I)}) \cong (0:_{\omega_{G(I)}} \m G(I)) =
(0:_{\omega_{G(I)}} \m ) = \oplus_{n \in \ZZ} (0:_{\omega_{n}}
\m).$$ (2) By definition $a(F(I))= -\min\{ n| [\omega_{F(I)}]_n \neq
0 \}$. Since $\omega_n$ is a finitely generated $A/I$-module for any
$n$ and $I$ is $\m$-primary, $A/I$ is of finite length and so 
$\omega_n$. As a consequence $(0:_{\omega_n}\m) \neq 0$ if and only
if $\omega_n \neq 0$. Therefore $a(F(I)) = a(G(I))$.

 On the other hand, it is known (see for instance (\cite{hz}, Proposition 3.6)) that if $G(I)$
is Cohen-Macaulay then $a(G(I))=r-d$, where $r:=r_J(I)$ is the reduction number of $I$ with
respect to any minimal reduction $J$ of $I$.

(3) Since $G(I)$ is Cohen-Macaulay, $G(I^k)$ is also Cohen-Macaulay for any positive integer
$k$ (see for instance (\cite{hz}, Corollary 4.6) ). On the other hand, by
(\cite{hz}, Corollary 4.6) $a(G(I^k))=[\frac{a(G(I))}{k}]$. Thus by (2)
$$a(F(I^k))= a(G(I^k))= [\frac{a(G(I))}{k}]= [\frac{a(F(I))}{k}]= [\frac{r-d}{k}].$$

(4) If $G(I)$ is Gorenstein then $\omega_{G(I)} \cong G(I)(a(G(I)))= G(I)(r-d)$. So
$\omega_n= I^{n+r-d}/I^{n+r-d+1}$ for any $n$ and by (1)
$$\omega_{F(I)} \cong \oplus_{n \in \ZZ} (0:_{\omega_n} \m) = \bigoplus_{n \in \ZZ} (I^{n+r-d+1}:\m) \cap I^{n+r-d}/I^{n+r-d+1}.$$
\end{proof}
 We know that if $G(I)$ and $F(I)$ are Cohen-Macaulay then $\reg G(I)=r=\reg F(I)$.
Now we prove that $\reg \omega_{G(I)}$ and $\reg \omega_{F(I)}$ are equal.
\begin{corollary}
Suppose $G(I)$ and $F(I)$ are Cohen-Macaulay. Then $\reg
\omega_{G(I)}= \reg \omega_{F(I)}$. In addition if $G(I)$ or $F(I)$
is Gorenstein then $\reg \omega_{G(I)} = d = \reg \omega_{F(I)}$.
\end{corollary}
\begin{proof}
Let $J$ be a minimal reduction of $I$ minimally generated
by $x_1,\ldots,x_d$ such that $x_1^{*},\ldots,x_d^{*} \in G(I)$ and
$x_1^{o},\ldots,x_d^{o} \in F(I)$ are regular sequences. Then $F(I)/J^o \cong F(I/J)$
and $G(I)/J^* \cong G(I/J)$.
Therefore $\reg \omega_{F(I)}= \reg (\omega_{F(I)}/J^0\omega_{F(I)})= \reg \omega_{F(I/J)}$
and $\reg \omega_{G(I)} = \reg \omega_{G(I/J)}$. But $\reg \omega_{F(I/J)}= a(\omega_{F(I/J)})$ and
$\reg \omega_{G(I/J)}= a(\omega_{G(I/J)})$. By Proposition \ref{pro4}(1) we have
$\omega_{F(I/J)}\cong \oplus_{n \in \ZZ}(0:_{[\omega_{G(I/J)}]_n}\m )$. From this we have
 $[\omega_{F(I/J)}]_n \neq 0$ if and only if $[\omega_{G(I/J)}]_n \neq 0$ for any $n$. Therefore
 $a(\omega_{F(I/J)}) = a(\omega_{G(I/J)})$. Thus $\reg \omega_{G(I)}= \reg \omega_{F(I)}$.

Now assume $G(I)$ is Gorenstein. Then $\omega_{G(I)} \cong
G(I)(r-d)$. Thus $\reg \omega_{G(I)}= \reg G(I)(r-d) = \reg
G(I)-r+d$. Since $G(I)$ is Cohen-Macaulay $\reg G(I)= r$. Therefore
$\reg \omega_{G(I)}= d$. If $F(I)$ is Gorenstein, then one can prove
the statement in a similar manner.
\end{proof}
\begin{corollary}
 Let $(A,\m)$ be a Gorenstein local ring of dimension $d > 0$
and $I$ is an $\m$-primary ideal with $G(I)$ is Cohen-Macaulay. Then
$\omega_{F(I)} \cong \bigoplus_{n \in \ZZ} \frac{(J^{n+r-d+1}:\m
I^r)\cap (J^{n+r-d}:I^r)}{(J^{n+r-d+1}:I^r)}$. If in addition $F(I)$
is Gorenstein, then $$\lambda\left( \frac{(J^{n+1}:\m I^r)\cap
(J^{n}:I^r)}{(J^{n+1}:I^r)}\right) = \lambda\left(\frac{I^n}{\m
I^n}\right)$$ for all
$n \in \ZZ^{+}$.
\end{corollary}
\begin{proof}
 By (\cite{hku}, Theorem 4.1), $\omega_{B}= \bigoplus_{n \in \ZZ} (J^{n+r}: I^r) t^{n+d-1}$, where
$B= A[It,t^{-1}]$ the extended Rees algebra of $I$. Since $\omega_{G(I)}= \omega_{B/t^{-1}B}= (\omega_B/t^{-1}\omega_B)(-1)$,
 we have $\omega_{G(I)}= \bigoplus_{n \in \ZZ} \frac{(J^{n+r-d}:I^r)}{(J^{n+r-d+1}:I^r)} $.
From Proposition \ref{pro4}(1),
\begin{eqnarray*}
 \omega_{F(I)} &=& \oplus_{n \in \ZZ} (0:_{[\omega_{G(I)}]_n} \m) \\
               &=& \bigoplus_{n \in \ZZ} \frac{(J^{n+r-d+1}:\m I^r)\cap
               (J^{n+r-d}:I^r)}{(J^{n+r-d+1}:I^r)}
\end{eqnarray*}
 Assume $F(I)$ is Gorenstein. Then $\omega_{F(I)} \cong F(I)(r-d)= \bigoplus_{n \in \ZZ} I^{n+r-d}/\m I^{n+r-d}$.
Therefore  $$\bigoplus_{n \in \ZZ} \frac{(J^{n+r-d+1}:\m I^r)\cap (J^{n+r-d}:I^r)}{(J^{n+r-d+1}:I^r)}=
\omega_{F(I)}= \bigoplus_{n \in \ZZ} I^{n+r-d}/\m I^{n+r-d}.$$ This implies that
$$\lambda\left( \frac{(J^{n+1}:\m I^r)\cap
(J^{n}:I^r)}{(J^{n+1}:I^r)}\right) =
\lambda\left(\frac{I^n}{\m I^n}\right)$$ for all $n \in \ZZ^{+}$.
\end{proof}
\begin{remark} \label{rmk1}
 Suppose $G(I)$ and $F(I)$ are Cohen-Macaulay rings. Let
$J=(x_1,\ldots,x_d)$ be a minimal reduction of $I$ such that
$x_1^*,\ldots,x_d^* \in G(I)$ and $x_1^o,\ldots,x_d^o \in F(I)$ are
superficial sequences and hence regular sequences. Denote
$J_i=(x_1,\ldots,x_i)$ and $J_0=(0)$. Then for any $i$ such that $1
\leq i \leq d$, $G(I)/(x_1^*,\ldots,x_i^*)\cong G(I/J_i)$ and
$F(I)/(x_1^o,\ldots,x_i^o) \cong F(I/J_i)$. Then $G(I)$ is Gorenstein if and
only if $G(I/J_i)$ is Gorenstein and
$F(I)$ is Gorenstein if and only if $F(I/J_i)$ is Gorenstein.
\end{remark}

Now we give a characterization for $F(I)$ to be Gorenstein in terms of
certain length conditions involving a minimal reduction of $I$ if $G(I)$
is Gorenstein.
\begin{theorem}\label{gor-char}
Let $(A,\m)$ be a Noetherian local ring, $I$ be an $\m$-primary ideal
and $J$ be a minimal reduction of $I$ with reduction number $r$.
Assume that $G(I)$ is a Gorenstein ring and $F(I)$ is
Cohen-Macaulay. Then $F(I)$ is Gorenstein  if and only if
$$\lambda\left(\frac{((I^{n+1}+J):\m)\cap
I^n}{I^{n+1}+JI^{n-1}}\right)=\lambda\left(\frac{I^n}{\m
I^n+JI^{n-1}}\right)$$ for $0 \leq n \leq r$.
\end{theorem}
\begin{proof}
Since $G(I)$ and $F(I)$ are Cohen-Macaulay, we may choose a generating
set for $J$ such that the corresponding images form a regular sequence
in $G(I)$ as well as in $F(I)$. Therefore we have
\begin{eqnarray*}
\omega_{F(I/J)} & = & (\omega_{F(I)}/J^o\omega_{F(I)})(r) \\
                &= &\bigoplus_{n \in \ZZ}\left[\frac{((I^{n+r+1}+J):\m)\cap
I^{n+r}+J}{I^{n+r+1}+J}\right].
\end{eqnarray*}
Using the isomorphism theorems and the fact that $G(I)$ is
Cohen-Macaulay, one obtains the isomorphism:
$$
[\omega_{F(I/J)}]_{n-r} \cong \frac{((I^{n+1}+J):\m)\cap
I^{n}}{I^{n+1}+JI^{n-1}}.
$$
Suppose $F(I)$ is Gorenstein. Then from the Remark \ref{rmk1}, it
follows that  $F(I/J)$ Gorenstein with canonical module
$\omega_{F(I/J)}=F(I/J)(r)$. Therefore
\begin{eqnarray*}
\lambda\left(\frac{((I^{n+1}+J):\m)\cap
I^n}{I^{n+1}+JI^{n-1}}\right)& =& \lambda([\omega_{F(I/J)}]_{n-r})\\
& = & \lambda([F(I/J)(r)]_{n-r})=\lambda\left(\frac{I^n}{\m
I^n+JI^{n-1}}\right)
\end{eqnarray*}
for all $n$.
Hence, in particular,
we get the required equality of lengths for
$0 \leq n \leq r$.
\vskip 2mm Conversely assume that $ \lambda(((I^{n+1}+J):\m)\cap
I^n/I^{n+1}+JI^{n-1})=\lambda(I^n/\m I^n+JI^{n-1})$ for $0 \leq n
\leq r$. Then by the above isomorphisms we have
$\lambda([\omega_{F(I/J)}]_n)=\lambda([F(I/J)(r)]_n)$ for all $n$.
This implies that
$\lambda(\omega_{F(I/J)})=\lambda(F(I/J)(r))=\lambda(F(I/J))$. Let $
\eta: P \rightarrow \omega_{F(I/J)}$ be the natural surjective map
from a graded free $F(I/J)$-module $P$ of rank equal to
$\mu(\omega_{F(I/J)})$, the minimal number of generators of
$\omega_{F(I/J)}$. Since $\omega_{F(I/J)}$ has finite injective
dimension, its injective dimension is equal to $\depth F(I/J) = 0.$
Therefore $\omega_{F(I/J)}$ is an injective module and hence
$\hom_{F(I/J)}(-,\omega_{F(I/J)})$ is an exact functor. Applying
this exact functor to the map $\eta,$ we get a surjective map
$\eta^*:\hom_{F(I/J)}(\omega_{F(I/J)},\omega_{F(I/J)})\longrightarrow
\hom_{F(I/J)}(P,\omega_{F(I/J)})$. But by the definition of
canonical module,
$$\hom_{F(I/J)}(\omega_{F(I/J)},\omega_{F(I/J)})\cong F(I/J).$$ Note
that $$\hom_{F(I/J)}(P,\omega_{F(I/J)})\cong \bigoplus_{rank(P)}
\omega_{F(I/J)}.$$ Hence there is a surjective map
$F(I/J)\rightarrow \bigoplus_{rank(P)} \omega_{F(I/J)}$. This
implies that
$$\lambda\left(\bigoplus_{rank(P)} \omega_{F(I/J)}\right) \leq
\lambda(F(I/J)).$$ This gives that $rank(P) \cdot
\lambda(\omega_{F(I/J)}) \leq \lambda(F(I/J))$. Since
$\lambda(\omega_{F(I/J)})=\lambda(F(I/J))$, the above inequality
gives that $rank(P)=1$. That is $\mu(\omega_{F(I/J)})=1$. Hence
$F(I/J)$ is Gorenstein with canonical module
$\omega_{F(I/J)}=F(I/J)(r)$. Since $J^o$ is generated by a regular sequence in
$F(I)$, $F(I)$ is Gorenstein with canonical module
$\omega_{F(I)} = F(I)(r-d)$.
\end{proof}

It is known that if $G(I)$ is Gorenstein, then $A$ is Gorenstein and
that such an analogue is not true in the case of fiber cone. We
show that if, in addition, $F(I)$ is Gorenstein, then $A/I$ is
Gorenstein.

\begin{corollary} \label{cor6}
Let $(A,\m)$ be a Noetherian local ring and $I$ is an $\m$-primary
ideal. Suppose $G(I)$ and $F(I)$ are Gorenstein. Then $A/I$
is Gorenstein.
\end{corollary}
\begin{proof}
 Put $n=0$ in the Theorem \ref{gor-char}, we get $\lambda((I:\m)/I)=\lambda(A/\m)=1$.
 This implies that $A/I$ is Gorenstein.
\end{proof}
Let $\pd_A(M)$ denote the projective dimension of $M$ as an $A$-module.
\begin{corollary} \label{cor7}
Suppose $I$ is an $\m$-primary ideal such that $G(I)$ and $F(I)$ are Gorenstein. 
Then $\pd(A/I) < \infty$ if and only if $I$ is generated by a regular
sequence.
\end{corollary}
\begin{proof}
Since $G(I)$ and $F(I)$ are Gorenstein, from the Corollary
\ref{cor6} it follows that $A/I$ is Gorenstein. Suppose $\pd(A/I)<
\infty$. Then from Theorem 2.6 of \cite{nv}, $I$ is generated by
a regular sequence. Conversely assume $I$ is generated by a regular
sequence. Then the Koszul complex $K_{\bullet}(I;A)$ is minimal free
resolution of $A/I$. Therefore $\pd(A/I)< \infty$.
\end{proof}
\begin{corollary}
  Suppose $(A,\m)$ is a regular local ring and $I$ is an $\m$-primary ideal of $A$ such that $G(I)$ and
$F(I)$ are Gorenstein. Then $I$ is generated by a regular sequence.
\end{corollary}
\begin{proof}
  By hypothesis $A$ is regular, so we have
$\pd(A/I) < \infty$. Therefore from Corollary \ref{cor7}, $I$ is
generated by a regular sequence.
\end{proof}

\begin{remark}
 Suppose $G(I)$ and $F(I)$ are Gorenstein rings and $\mu(I)= d+2$. Then by Corollary
 \ref{cor6}, $A/I$ is Gorenstein. By (\cite{hku}, Remark 2.9(1),
 (2)) we have $\pd_A(G(I)) < \infty$ if and only if $G(I)$ is a
 complete intersection. That is the defining ideal of $G(I)$ is
 generated by a regular sequence. In this case $F(I)$ is also a
 complete intersection.
\end{remark}

The expression we have obtained for the canonical module of the fiber
cone does not reveal when it can be realized as a submodule of
$F(I)$. Below we give a sufficient condition for $\omega_{F(I)}$
to be a submodule of $F(I)$.

\begin{corollary}
Let $(A, \m)$ be a Noetherian local ring and $I$ be an $\m$-primary
ideal. Assume that $G(I)$ is a
Gorenstein ring. Suppose $F(I)$ is Cohen-Macaulay and $A/I$ is Gorenstein.
If $(I^{n+r-d+1}:\m) \cap \m I^{n+r-d}= I^{n+r-d}$ for all $n \geq d-r+1$,
then $\omega_{F(I)}$ is an ideal of $F(I)$.
\end{corollary}
\begin{proof}
For any $n \geq d-r$, there is a natural map
$\psi_n:(I^{n+r-d+1}:\m) \cap I^{n+r-d}/I^{n+r-d}\longrightarrow
I^{n+r-d}/\m I^{n+r-d}$ which gives rise to a natural $F(I)$-linear
map $\psi : \omega_{F(I)} \longrightarrow F(I)$. The kernel of
$\psi_n$ is $(I^{n+r-d+1}:\m) \cap \m I^{n+r-d}/I^{n+r-d}$. Then the
$F(I)$-linear map $\psi:\omega_{F(I)}\longrightarrow F(I)(r-d)$ has
kernel $\oplus_{n \geq d-r}(I^{n+r-d+1}:\m) \cap \m
I^{n+r-d}/I^{n+r-d}$. Then $\omega_{F(I)}$ is an ideal of $F(I)$  if
$ker(\psi)=0$, i.e., if $(I^{n+r-d+1}:\m) \cap \m I^{n+r-d}=
I^{n+r-d}$ for all $n \geq d-r+1$. This completes the proof.
\end{proof}
The following Proposition shows that $e_0(\omega_{F(I)}) <
e_0(\omega_{G(I)})$ unless $I = \m$.
\begin{proposition}
Let $(A, \m)$ be a Noetherian local ring and $I$ be an $\m$-primary
ideal. Suppose $G(I)$ and $F(I)$ are Cohen-Macaulay. Then
$e_0(\omega_{F(I)}) \leq e_0(\omega_{G(I)})$. Furthermore, if $G(I)$
is Gorenstein, then equality holds if and only if $I=\m$.
\end{proposition}
\begin{proof}
From Proposition \ref{pro4}(1), it is clear that
$\lambda([\omega_{F(I)}]_n) \leq \lambda([\omega_{G(I)}]_n)$ for all
$n$. Hence $e_0(\omega_{F(I)}) \leq e_0(\omega_{G(I)})$. Let
$J=(x_1,\ldots,x_d)$ be a minimal reduction of $I$ with reduction
number $r$ such that $x_1^o,\ldots,x_d^o \in F(I)$ is a regular
sequence for $F(I)$, $\omega_{F(I)}$ and $x_1^*,\ldots,x_d^* \in
G(I)$ is $G(I)$ and $\omega_{G(I)}$ regular sequence. Then
$e_0(\omega_{G(I)})=e_0(\omega_{G(I/J)})=\lambda(\omega_{G(I/J)})$
and
$e_0(\omega_{F(I)})=e_0(\omega_{F(I/J)})=\lambda(\omega_{F(I/J)})$.

 Now assume $e_0(\omega_{F(I)}) = e_0(\omega_{G(I)})$. Then
 $\lambda(\omega_{F(I/J)}) = \lambda(\omega_{G(I/J)})$. By
 Proposition \ref{pro4}(1) we have $\omega_{F(I/J)} \subseteq
 \omega_{G(I/J)}$. Therefore $\omega_{F(I/J)} =\omega_{G(I/J)}$.
 This implies that $(I:\m)/I=A/I$. This gives
$\lambda((I:\m)/I)= \lambda(A/I)$. This implies $(I:\m)=A$. That is
$I=\m$. The converse always holds. This completes the proof.
\end{proof}

We conclude this article by providing two examples. First we give an
example of an ideal whose associated graded ring is Gorenstein but the
fiber cone is not.
\begin{example}
Let $A=k[\![t^4,t^9,t^{10}]\!]$ and $I=(t^8,t^{18},t^{10}),\m =
(t^4,t^9,t^{10})$. Then $A$ is a one dimensional Gorenstein local
domain and $I$ is an $\m$-primary ideal. $J=(t^8)$ is a minimal
reduction of $I$ of reduction number $1$. Then it follows from the
Corollary 4.5 (5) of \cite{hku} that $G(I)$ is Gorenstein. Since $I$
has reduction number $1$, $F(I)$ is Cohen-Macaulay. Since $\mu(I) = 3
> \dim A + 1$, by Proposition 4.1 of \cite{jpv}, $F(I)$ is not
Gorenstein. Also, note that
$$
5 = \lambda\left (\frac{((I^2+J):\m) \cap I}{I^2+J} \right ) \neq
\lambda \left ( \frac{I}{\m I+J} \right )= 1.
$$
\end{example}
In the example below we apply our result to obtain an example of a
Gorenstein fiber cone.
\begin{example}[\cite{hku}, Example 2.5]
Let $A=k[\![t^4,t^9,t^{10}]\!]$ and $I := (t^8,t^9,t^{10})$.
Then $A$ is a one dimensional Gorenstein local
domain, $I$ is an $\m$-primary ideal and $J=(t^8)$ is a minimal
reduction of $I$ with reduction number $2$. Let $R=A/I[X,Y,Z]$. Then
$G(I) \cong \frac{R}{(XZ-Y^2, wX^2-Z^2)}$, where $w$ is the image of
$t^4$ in $A/I$. Since $XZ-Y^2, wX^2-Z^2$ is an $R$-regular sequence
and $A/I$ is Gorenstein therefore $G(I)$ is Gorenstein. Since $\m I^n
\cap J = \m JI^{n-1}$ for all $n \geq 1$,  $F(I)$ is Cohen-Macaulay.
By using CoCoA, \cite{co}, it can easily be seen that
$$\lambda\left (\frac{((I^2+J):\m) \cap I}{I^2+J} \right )=
2 = \lambda \left ( \frac{I}{\m I+J} \right ) $$ and $$\lambda\left
(\frac{((I^3+J):\m) \cap I^2}{I^2+JI} \right )= 1= \lambda \left (
\frac{I^2}{\m I^2+JI} \right ).$$ Therefore by the Theorem
\ref{gor-char}, $F(I)$ is Gorenstein.
\end{example}


\begin{thebibliography}{AAAA}

\bibitem [BH]{bh} W. Bruns, J. Herzog, {\em Cohen-Macaulay
rings}, Revised Edition, Cambridge Studies in Advanced Mathematics,
39. Cambridge University Press, Cambridge, 1998.

\bibitem [CH]{ch} A. Conca, J. Herzog, {\em Castelnuovo-Mumford regularity
of products of ideals}, Collect. Math. {\bf 54} (2003), 137-152.

\bibitem[C]{co2} T. Cortadellas, {\em Fiber cones with almost maximal depth},  Comm. Algebra  {\bf 33}  (2005),  no. 3, 953--963.

\bibitem[CZ1]{cz2} T. Cortadellas, S. Zarzuela, {\em On the depth of the fiber cone of filtrations},  J. Algebra  {\bf 198}  (1997),  no. 2, 428--445.

\bibitem[CZ2]{cz1} T. Cortadellas, S. Zarzuela, {\em Depth formulas for certain graded modules associated to a filtration: a survey}, Geometric and combinatorial aspects of commutative algebra (Messina, 1999), 145--157, Lecture Notes in Pure and Appl. Math., 217, Dekker, New York, 2001.

\bibitem[CZ3]{cz3} T. Cortadellas, S, Zarzuela, {\em On the structure
of the fiber cone of ideals with analytic spread one}, J. Algebra {\bf
317} (2007), no. 2,  759--785.

\bibitem[Co]{co} CoCoATeam, {\em CoCoA: a system for doing
  Computations in Commutative Algebra}. Available at
  http://cocoa.dima.unige.it.

\bibitem[GI]{gi} S. Goto and S.-I. Iai, {\em Embeddings of certain
  graded rings into their canonical modules}, J. Algebra {\bf 228}
  (2000), 377--396.

\bibitem[GN]{gn} S. Goto, K. Nishida, {\em The Cohen-Macaulay and Gorenstein Rees algebras associated to filtrations}, Mem. Amer. Math. Soc. {\bf 110} (1994), no. 526. American Mathematical Society, Providence, RI,  1994.

\bibitem[HIO]{hio} M. Herrmann, S. Ikeda, and U. Orbanz, {\em Equimultiplicity and blowing up}, Springer-Verlag, Berlin, 1988.

\bibitem[HZ]{hz} L. T. Hoa and S. Zarzuela, {\em Reduction number and
  a-invariant of good filtrations}, Comm. Algebra {\bf 14} (1994), no. 22, 5635--5656.

\bibitem [Hu]{h} C. Huneke, {\em The theory of d-sequences and powers of ideals}, Advances in Math. {\bf 46}
 (1982), 249--279.

\bibitem [Hy]{eh} E. Hyry, {\em On the Gorenstein property of the associated
graded ring of a power of an ideal}, Manuscripta Math. {\bf 80}
(1993), 13--20.

\bibitem [HKU]{hku} W. Heinzer, M.-K. Kim, B. Ulrich, {\em The
  Gorenstein and complete intersection  properties of associated
  graded rings}, J. Pure and Applied Algebra {\bf 201} (2005), 264--283.

\bibitem[HRS]{hrs} M. Herrmann, J. Ribbe, P. Schenzel, {\em On the
  Gorenstein property of the form rings}, Math. Z. {\bf 213} (1993),  301--309.

\bibitem [JPV]{jpv} A. V. Jayanthan, Tony J. Puthenpurakal, J. K.
  Verma, {\em On fiber cones of $\m$-primary ideals}, Canad. J. Math.
  {\bf 59} (2007), no.1, 109--126.

\bibitem[NV]{nv} S. Noh and W. Vasconcelos, {\em The $S_2$-closure of
  a Rees algebra}, Results. Math. {\bf 23} (1993), 149 -- 162.

\bibitem [O]{o} A. Ooishi, {\em Genera and arithmetic genera of commutative rings}, Hiroshima Math. J. {\bf 17} (1987), 47--66.

\bibitem [T]{t1} N. V. Trung, {\em The Castelnuovo regularity of the
Rees algebra and the associated graded ring},  Trans. of Amer. Math.
Soc. {\bf 350} (1998), no. 7,  2813--2832.

\bibitem[TVZ]{tvz} N. V. Trung, D. Q. Vi\^{e}t, S. Zarzuela, {\em When is the Rees algebra Gorenstein?}, J. Algebra {\bf 175} (1995), 137--156.
\end{thebibliography}
\end{document}